\documentclass[12pt,a4paper]{amsart}
\usepackage{lineno}
\usepackage{amssymb}
\usepackage[final]{showkeys}
\usepackage{microtype}
\usepackage{color}
\usepackage{graphicx}
\usepackage[hmargin=3.5cm,vmargin={5cm,5cm}]{geometry}
\numberwithin{equation}{section}

\allowdisplaybreaks

\begin{document}

    \title[]{A note on a modified fractional Maxwell model}
    
    \author{R. Garra$^1$}
    \address{${}1$ Institute of Marine Sciences, National Research Council (CNR), Rome, Italy\\
    e-mail: rolinipame@yahoo.it}

	\author{A. Consiglio$^2$}
	\address{$^2$ Institut f\"{u}r Theoretische Physik und Astrophysik and W\"{u}rzburg-Dresden Cluster of Excellence ct.qmat, Universit\"{a}t W\"{u}rzburg, 97074 W\"{u}rzburg, Germany}
    
    \author{F. Mainardi$^3$}
    \address{${}^3$ Dipartimento di Fisica e Astronomia, Università di Bologna, \& INFN,
    	Via Irnerio 46, I-40126 Bologna, Italy.}

    \keywords{Linear viscoelasticity, creep, relaxation,  Hadamard fractional derivative, fractional calculus, ultra slow kinetics.
     \\
    {\it MSC 2010\/}:  
	26A33, %% Fractional Derivatives and integrals
	%44A10, %% Laplace Transforms 
    45D05, %% 	Voltera Integral Equations
	74D05, %% Linear consitutive equations
	74L10, %% Soil and Rock Mechanics
	76A10  %5% Viscoelastic fluids
    }

    \date{\today}

    \begin{abstract}
     
    In this paper we consider a modified fractional Maxwell model based on the application of Hadamard-type fractional derivatives.
    The model is physically motivated by the fact that we can 
    take into account at the same time memory effects and the time-dependence of the viscosity coefficient. 
    We obtain an ultra-slow relaxation response whose explicit analytic form is given by the Mittag-Leffler function with a logarithmic argument. We show graphically the main properties of this relaxation response, also with the asymptotic behaviour.

	\smallskip
	\noindent
{\bf 	This paper has been published in Chaos, Solitons and Fractals, Vol 103 (2022) 112544/1-5. 
	DOI: 10.1016/j.chaos.2022.112544}
	
    \end{abstract}

    \maketitle

    \section{Introduction}
    
	Classical viscoelastic models are based on the combination of viscous dashpots and springs (in series or parallel), leading to different 
	governing equations. The viscosity is considered as a constant in the classical models but in recent papers a time-varying viscosity $\eta(t)$ was
	introduced, we refer in particular to Pandey and  Holm \cite{Holm} and Yang et al.\cite{Yang}. 
	According to \cite{Yang} a modified Maxwell model based on a time-varying viscosity can be formulated as
	\begin{equation}
		\sigma+\frac{\eta(t)}{E}\frac{d\sigma}{dt} = \eta(t)\frac{d\epsilon}{dt},
	\end{equation}
    where $\sigma$ and $\epsilon$ are the stress and the strain respectively, while $E$ is the elastic modulus of the spring.\\
   The authors in \cite{Holm}, \cite{Yang}   have considered, in particular, a linearly time-varying viscosity 
    $$\eta(t) = \eta_0 +\theta t,$$
    where 
    $\eta_0$ denotes the initial viscosity and $\theta$ the strain-hardening coefficient. \\
    In this case the modified Maxwell model becomes
    \begin{equation}\label{mod}
 	\sigma+\frac{\eta_0 +\theta t}{E}\frac{d\sigma}{dt} = \left(\eta_0 +\theta  t\right)\frac{d\epsilon}{dt}.
    \end{equation}
    On the other hand, it is well-known that viscoelastic and mechanical models with memory based on the applications of fractional derivatives and integrals are widely used and also validated (we refer for example to \cite{spada1}, \cite{ps} and the references therein). \\
    In the recent paper by Garra, Mainardi and Spada \cite{spada},  the authors 
    have considered a modified Lomnitz model based on the application of Hadamard-type integro-differential operators.
   % \textcolor{red}
   {Here we consider a little modification of the operator considered in \cite{spada} that is, for $\nu \in(0,1)$
    \begin{equation}
    	\begin{array}{ll}
    		\widehat{O}^t_\nu f(t) &:={\displaystyle
    			\frac{1}{\Gamma(1-\nu)}} \left(\frac{E}{\theta}\right)^{-\nu}\displaystyle
    			\int_{\frac{1-\eta_0}{\theta}}^{t}	
    			\!\!  
    			\ln^{-\nu}\left(\frac{\eta_0+\theta t}{\eta_0+\theta\tau}\right)
    			 \left(\frac{\eta_0 +\theta \tau}{E}\frac{d}{d\tau}f(\tau)\right)\frac{E}{\eta_0+\theta t}\, d\tau, 
    	\end{array}
    \end{equation}
    that provides a fractional generalization of the differential operator appearing in the governing equation of the modified Maxwell model equation considered in \cite{Yang}. 
    Indeed, for $\nu = 1$
    we recover the differential operator
    \begin{equation}
    \left(\frac{\eta_0 +\theta t}{E}\frac{d}{dt}\right)
\end{equation}
      appearing in the model equation \eqref{mod}.
  This is a special case of a fractional derivative w.r.t. another function obtained by means of the deterministic change of variable
$t\rightarrow \frac{E}{\theta}\ln\left(\eta_0+\theta t\right)$ (see \cite{almeida}).}
       We observe that, by taking such deterministic time-change, we can consider the time-dependence of the viscosity in the classical model. By using the Hadamard-type derivative, we have this change of variable inside the Caputo fractional derivative. 
    We take inspiration from the analysis developed in \cite{spada} to consider a completly new generalization of the fractional Maxwell model. Indeed, here we consider the relaxation response of this new theoretical model that is not coincident with the Lomnitz-type profile obtained in \cite{spada}.
    The motivation of this mathematical manipulation is clearly motivated by the fact that fractional models play a relevant role in the field of viscoelasticity (see e.g. \cite{Mainardi_BOOK10} and the references therein). \\
    In this paper we introduce a fractional modified Maxwell model based on Hadamard-type derivatives and we find the explicit form of the relaxation response. Then, we study the mathematical properties of the obtained result, providing a short discussion about the asymptotic behaviour and the graphical representation by using the recent numerical routine developed by 
  Garrappa in \cite{Garrappa MATLAB14}.  As far as we know the most efficient and accurate routine 
for the computation of the Mittag-Leffler function (in the whole complex plane)
 is indeed that developed by Garrappa in 
\cite{Garrappa MATLAB14} and illustrated in \cite{Garrappa SIAM15} published by \textit{SIAM Journal of Numerical Analysis}. 
%% and based on the numerical inversion of the Laplace transform.
% This routine allows us  to compute the Mittag-Leffler function with up to 3 independent %parameters (the so-called Prabhakar function).
    
    We finally show that the modified Maxwell model considered in \cite{Yang} can be recovered as a special case.
    
    \section{Preliminaries on Hadamard-type fractional derivatives \label{ope}}
    
   In the paper 
   by Beghin, Garra and Macci \cite{jap}, an integro-differential operator with logarithmic kernel has been introduced in the context of correlated fractional negative binomial processes in statistics. 
   Then, the application of these operators in relation to the generalized Lomnitz law was discussed in \cite{spada}.
   % \textcolor{red}
   {Here we consider a little modification of the operator considered in \cite{spada} motivated by the model that we are considering. The time-evolution operator 
    $\widehat{O}^t_\nu$  acting on a sufficiently well-behaved function 
    $f(t)$, for $\nu \in (0,1)$
    is defined as
\begin{equation}
	\begin{array}{ll}
		\widehat{O}^t_\nu f(t) &:={\displaystyle
			\frac{1}{\Gamma(1-\nu)}} \left(\frac{E}{\theta}\right)^{-\nu}\displaystyle
		\int_{\frac{1-\eta_0}{\theta}}^{t}	
		\!\!  
		\ln^{-\nu}\left(\frac{\eta_0+\theta t}{\eta_0+\theta\tau}\right)
		\left(\frac{\eta_0 +\theta \tau}{E}\frac{d}{d\tau}f(\tau)\right)\frac{E}{\eta_0+\theta t}\, d\tau.
	\end{array}\label{defO}
\end{equation}    
    We can consider this operator as the fractional power of the operator appearing in the modified Maxwell model, namely
    \begin{equation}
    \left(\frac{\eta_0 +\theta t}{E}\frac{d}{dt}\right)^\nu f(t) = 	\widehat{O}^t_\nu f(t).
    \end{equation}
    }
    Moreover, for $\nu = 1$ (as usual in the context of fractional differential operators), we recover the ordinary differential operator with a variable coefficient
       \begin{equation}
    \widehat{O}_{\nu= 1}^t f(t) =	\left(\frac{\eta_0 +\theta t}{E}\frac{d}{dt}\right) f(t).
    \end{equation}
    A relevant property of this operator is given by the following result 
    \begin{equation}\label{pr}
    \widehat{O}^{t}_\nu \ln^\beta (\eta_0+\theta t) =\left(\frac{E}{\theta}\right)^{-\nu} \frac{\Gamma(\beta+1)}{\Gamma(\beta+1-\nu)}\ln^{\beta-\nu}(\eta_0+\theta t)
    \end{equation}
    for $\nu \in (0,1)$ and $\beta>-1 \setminus{\{0\}}$,
      see \cite{almeida}, for the details.\\
    Moreover we have that 
    \begin{equation}
     \widehat{O}^{t}_\nu \ const. = 0.
     \end{equation}
%%%

    In the general framework of the classical theory of  fractional calculus (see e.g. the monograph by Kilbas, Srivastava and Trujillo \cite{kilbas}), we observe that these Hadamard-type fractional derivatives can be considered as a special case of fractional derivative of a function with respect to another function, recently named in the literature $\psi$-fractional derivatives (we refer for example to the recent papers \cite{almeida} and \cite{alm}).

%     \textcolor{red}
{Indeed, we can observe that, heuristically, the integro-differential operator 
    $\widehat{O}_\nu^t$ can be obtained by means of a 
    time-change $t\rightarrow \frac{E}{\theta}\ln\left(\eta_0+\theta t\right)$ starting from the definition of the Caputo fractional derivative.}
    We also note  that the operator $\widehat{O}^t_{\nu}$ can be considered as a sort of \textit{fractional} counterpart of the differential operator $\widehat{O}_1^t$ appearing in the stress-strain equation governing the modified Maxwell equation with time-varying viscosity. 
    We remark that for $\eta_0= 0$ and $\theta=1$ 
    the operator $\widehat{O}^t_\nu$ coincides with the regularized Hadamard fractional derivative (see \cite{jap} for more details).
    
    We finally observe that, by using the property \eqref{pr} (see also \cite{almeida}), it can be  shown that the composed Mittag-Leffler function 
    \begin{equation}\label{func}
    E_{\nu,1}\bigg(-\left(\frac{E}{\theta}\right)^\nu\ln^\nu(\eta_0+\theta t)\bigg)= \sum_{k=0}^\infty\frac{\bigg(-\displaystyle\left(\frac{E}{\theta}\right)^\nu\ln^\nu(\eta_0+\theta t)\bigg)^k}{\Gamma(\nu k+1)}
    \end{equation}
    is an eigenfunction of the operator $\widehat{O}^t_\nu$.
Once again we find  the   Mittag-Leffler functions, that are known to  play a central role in the theory of fractional differential equations, see  e.g.  the recent monograph by Gorenflo et al. \cite{book}. 
%%%%%
As a matter of fact, the composition of the classical Mittag-Leffler function with the power of a logarithmic function decays as a power of the logarithmic function.
 This behaviour can be interesting 
in the framework the so-called {\it ultra-slow kinetics} that includes 
phenomena of strong anomalous relaxation and diffusion and related stochastic processes.
On this respect, we refer for example to 
Chechkin et al. \cite{Chechkin-et-al_EPL03}, 
Metzler and Klafter \cite{Metzler-Klafter_JPhysics04},
Mainardi et al. \cite{Mainardi-et-al_JVC07,Mainardi-et-al_JVC08}  
  and more recently to Wen Chen et al. \cite{WenChen-et-al_FCAA16}.
Ultra-slow decay in luminescence have been recently studied in \cite{gp}
leading to the function \eqref{func}.

    \section{The fractional modified Maxwell model}
    
    We here consider the following fractional modified Maxwell model 
  	\begin{equation}
  	\sigma+ \widehat{O}^{t}_\nu \sigma = \widehat{O}^{t}_\nu \epsilon,
  	\end{equation}
  with $\nu \in (0,1)$. For $\nu = 1$ we will show that we recover the modified Maxwell model considered in \cite{Yang} (with a suitable choice of the parameters appearing in the operator).\\
  Let us consider, in particular, the relaxation test, where the strain is set to be a constant $\epsilon = \epsilon_0$, for $t>0$, then
  $\widehat{O}^{t}_\nu \epsilon = 0$ and we have the relaxation equation
  	\begin{equation}
  \widehat{O}^{t}_\nu \sigma = -\sigma.
  \end{equation}
  Then, we have that the relaxation response is given by (see the previous section)
  \begin{equation}\label{sigma-Hadamard}
\sigma(t) = E_{\nu,1}(-\ln^\nu(1+t)).
  \end{equation}
In the special case $\nu = 1$, we have
\begin{equation}
	\sigma(t) = \exp\{-\ln(1+t)\} = (1+t)^{-1},
\end{equation}
that coincides with the relaxation response in \cite{Yang}, eq.(8), where we have taken for simplicity $E = \theta = \eta_0 = 1$.

% \textcolor{red}
{In Figure 1 and Figure 2} we present some illuminating plots of the composite Mittag-Leffler function 
 \eqref{sigma-Hadamard}
 corresponding to   rational values of the parameter $\nu$ 
 ($\nu = 0.25, 0.50, 0.75, 0.9, 1$) by adopting
 four different scale. This is for possible convenience of the experimental researchers that measure the  stress relaxation in viscoelasticity.
As outlined in the introduction we use the Matlab routine of Garrappa 
\cite{Garrappa MATLAB14}.  

 \begin{figure}  %%[H]
\centering
\includegraphics[width=0.495\textwidth, trim={28cm  1cm  28cm 2.8cm},clip]{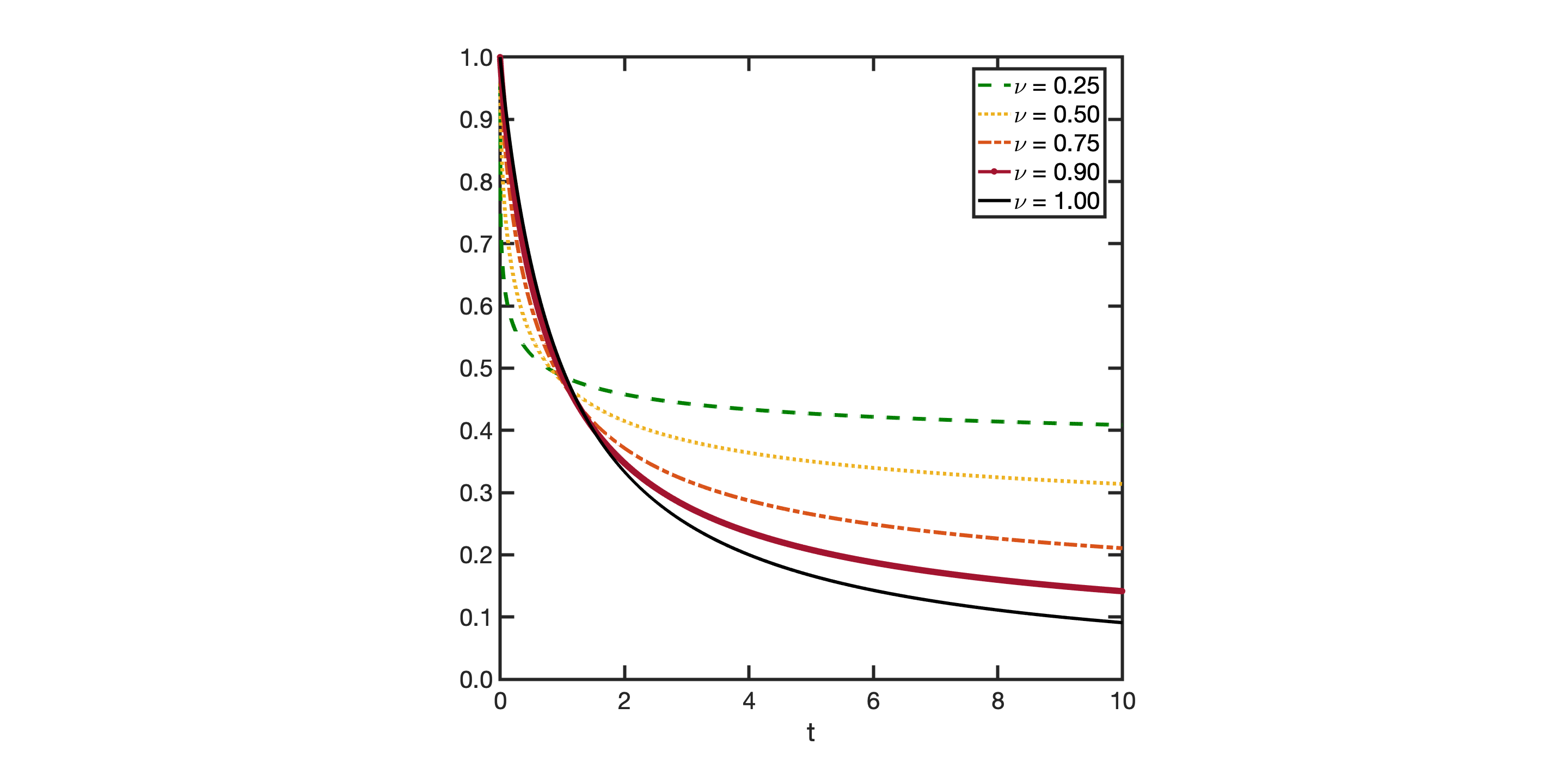}
\includegraphics[width=0.495\textwidth, trim={28cm  1cm  28cm 2.8cm},clip]{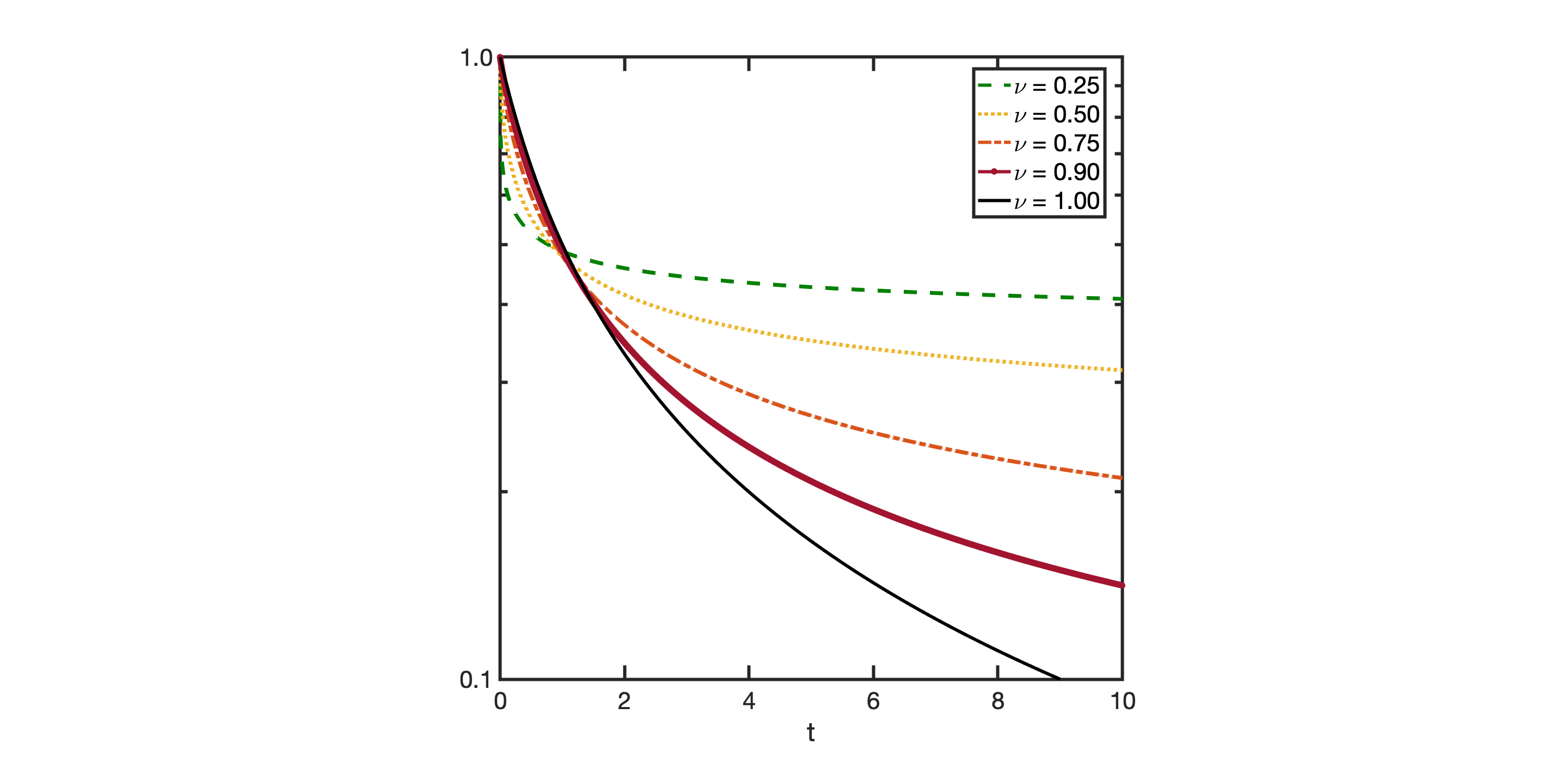}
\caption{The composed  function 
$E_{\nu,1}(-\ln^\nu(1+t))$
for $\nu=0.25, 0.50, 0.75, 0.9, 1$ in the time range  $0\le t \le 10$.}
\label{fig1}
\end{figure}
 \begin{figure}  %%[H]
\centering
\includegraphics[width=0.495\textwidth, trim={28cm  1cm  28cm 2.8cm},clip]{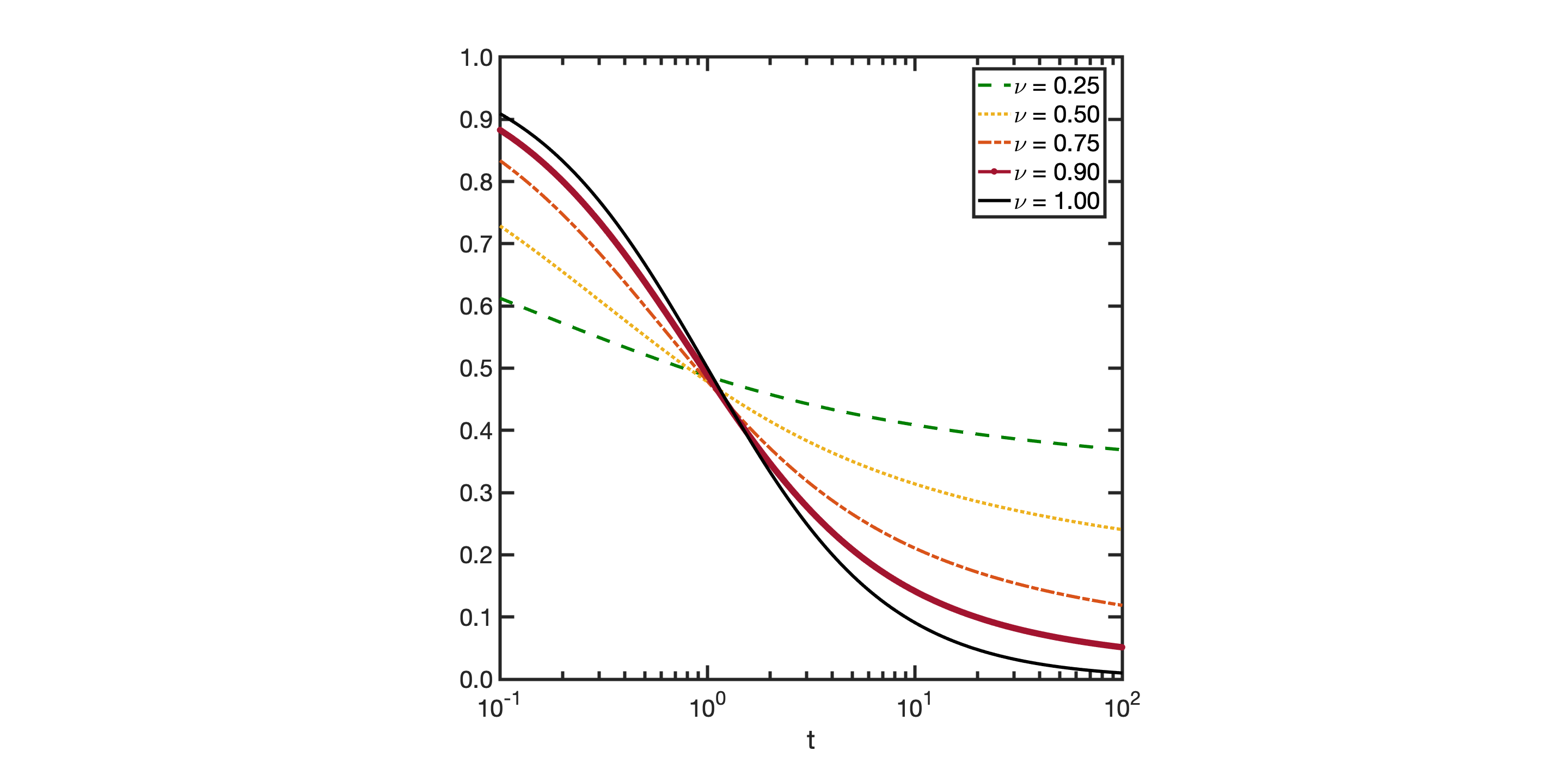}
\includegraphics[width=0.495\textwidth, trim={28cm  1cm  28cm 2.8cm},clip]{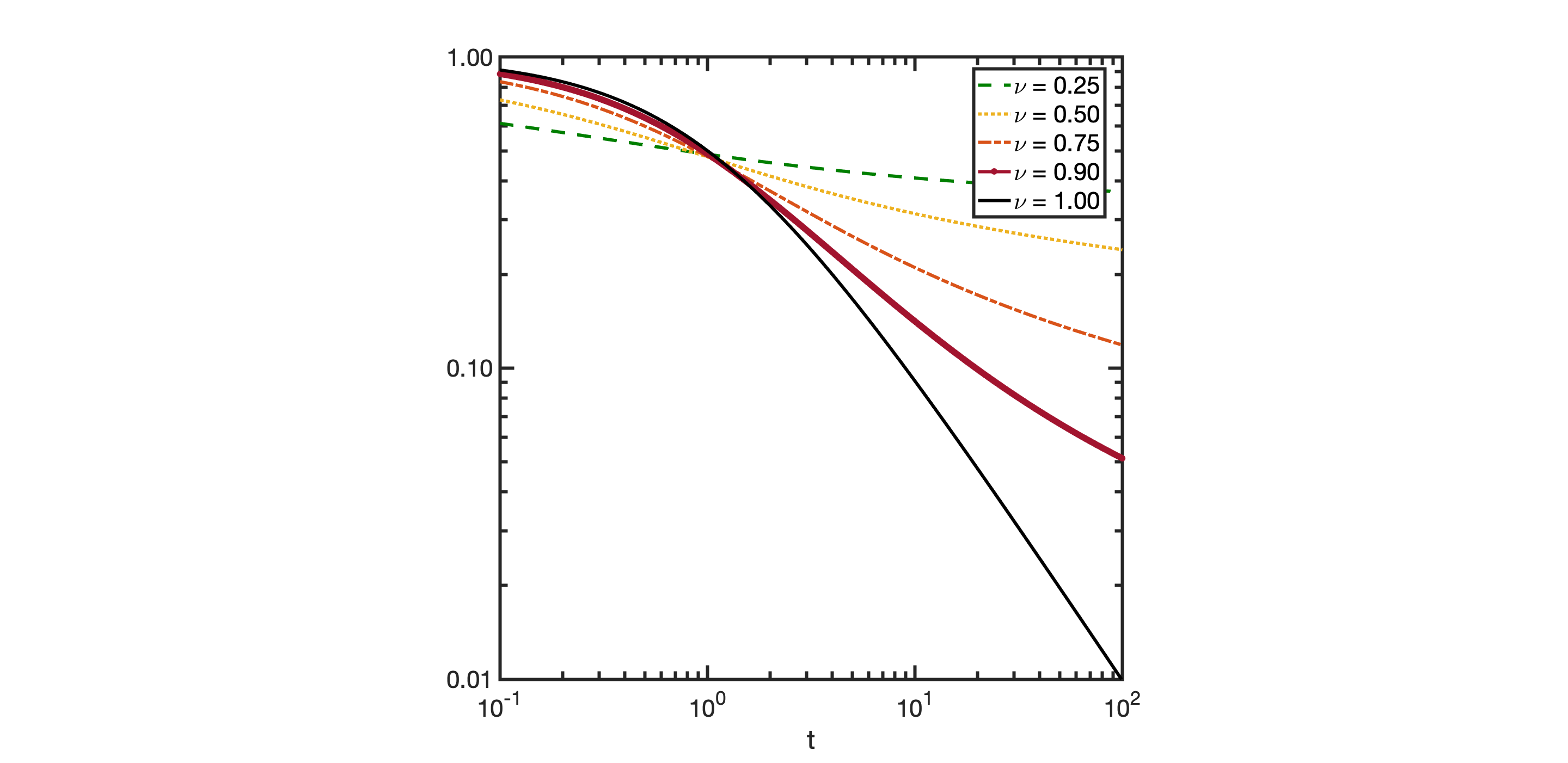}
\caption{The composed  function 
$E_{\nu,1}(-\ln^\nu(1+t))$
for $\nu=0.25, 0.50, 0.75, 0.9, 1$ in the time range  $10^{-1}\le t \le 10^{+2}$.}
\label{fig2}
\end{figure}

The asymptotic representation for $t\to\infty$ of the function 
$E_{\nu,1}(-\ln^\nu(1+t))$ is readily obtained by the first term of the asymptotic
series of the Mittag-Leffler function for $x>0$ and $\nu \in (0,1)$, see e.g \cite{book}:
\begin{equation}
E_{\nu, 1}(-x) \sim  -\sum_{n=1}^\infty\dfrac{(-1)^{n-1}\, x^{-n}}{\Gamma (1-\nu n)},
\quad  x\to\infty.
\end{equation}
Then 
\begin{equation}
 E_{\nu,1}(-\ln^\nu (1+t)) \sim \dfrac{\ln^{-\nu}(1+t)}{\Gamma(1-\nu)}, \quad
 t\to\infty.
 \end {equation}
 so that
 we also try to find the matching with the above asymptotic representations as $t \to \infty$.
 We note that this matching is expected to be for very large values of time because
 the presence of the logarithm that is known to be a {\it slow varying function}.
 Of course in the limit  $\nu \to 1^-$ we get a singular transition  from a logarithmic decay to a 
 power law decay because 
\begin{equation}
 E_{1,1}(-\ln (1+t)) =\exp(-\ln (1+t)) = \dfrac{1}{1+t}\,, \quad t>0\,.
 \end{equation}
 Then for this matching we add the case $\nu=0.95$ to point out the singular limit
 as $\nu$ is closer to 1. 
  \begin{figure}  %%[H]
 	\centering
 	\includegraphics[width=0.99\textwidth, trim={0.06cm  1cm  0.06cm 0.8cm},clip]{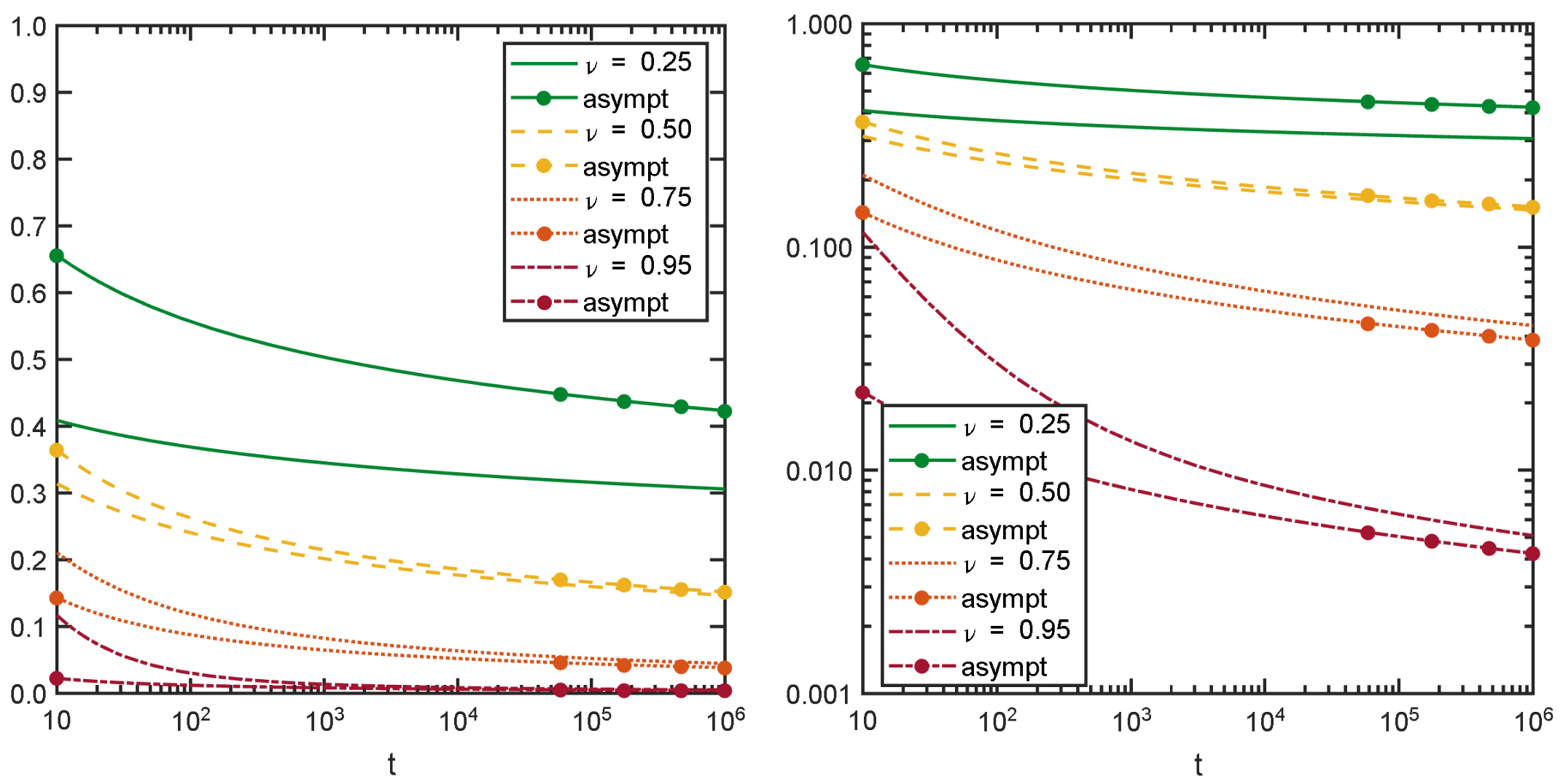}
 	\caption{Plots in an enlarged logarithmic scale, compared with their asymptotic representations.}
 	\label{fig3}
 \end{figure}
 This can be understood in analogy with the singular limit of the classical
 Mittag-Leffler function $E_\nu(-x)$  of negative argument when $\nu \to 1^-$ 
 as described
 in the recent paper by Paris \cite{Paris FCAA22} published in \textit{Fractional Calculus and Applied Analysis}
 to which the interested reader is asked to refer. %% follow the analysis in \cite{Paris FCAA22}.
 %\textcolor{red}
 {In Figure 3, we show the comparison between the plots in an enlarged logarithmic scale and their asymptotic representations. It is evident that the accurate value of the Mittag-Leffler function approaches the corresponding asymptotic representation with a good agreement.}
\\ 
Since the function $\ln^\nu(1+t)$  for $\nu \in(0,1)$ is a Bernstein function
 as in the case $\nu=1$, then the stress 
\eqref{sigma-Hadamard} 
 is relaxing as a completely monotone function in analogy  with standard fractional Maxwell model,
 see e.g.\cite{Mainardi_BOOK10}. Indeed, we observe that this function can be obtained as the composition of the two Bernstein functions $x^\nu$ and $log (1+t)$, see e.g. the treatise by Schilling et al \cite{sch}.

%%%% 

   \section{Conclusions}
   
   In this paper we consider an explorative mathematical modification of the fractional Maxwell model based on the application of Hadamard-type derivatives. It is well-known and widely accepted that time-fractional models in viscoelasticity have the advantage to take into account memory effects in 
   relaxation processes. The main aim of this paper is to combine two different effects in a single theoretical model. Considering the recent studies about the role of time-dependent viscosity coefficient in the Maxwell model (see \cite{Holm}), we here suggest to consider in a single operator memory-effects and time-dependent viscosity. 
   This fractional operator results as a deterministic time-change from the Caputo derivative and the modified Maxwell model leads to an ultra-slow relaxation process. 
   We underline that the relaxation response obtained in \cite{Holm} can be recovered as a special case. \\
   The idea to apply Hadamard-type derivatives in viscoelasticity is inspired by the recent studies for the formulation of a generalized Lomnitz creep law (see \cite{spada}). Moreover, different mathematical approaches for ultra-slow relaxation models are gaining interest and this topic is of relevant interest for the applications (we refer in particular to
   % \textcolor{red}
   {\cite{magin} and \cite{rec}}).
   In our view the analysis of ultra-slow relaxation models based on the Hadamard-fractional calculus approach  should be still deepened. 
   From the theoretical point of view, Hadamard-type derivatives can be obtained as a particular choice of fractional derivatives w.r.t. another functions. Modified Maxwell models based on fractional derivatives w.r.t. another functions can be object of future interest in order to consider different time-dependence of the viscosity coefficients and leading to different relaxation response.
   A similar idea was recently considered for generalized Scott-Blair models in \cite{andre}.\\
   We finally observe that a similar model was considered in \cite{long}. In comparison with this paper, we introduce the analysis of this model in the context of the literature, providing physical motivations and new graphical simulations.
  
   \section*{Acknowledgments}
   We are greatful to both referees for their comments and suggestions.
	The work of the authors %% R. G., A. C. and F. M.
	 has been carried out in the framework of the activities of the National Group of Mathematical Physics (GNFM, INdAM).
	The authors are grateful to Prof R.B Paris and to Prof. T. Simon for helpful comments. 
	FM would like to thank the Isaac Newton Institute for
	Mathematical Sciences, Cambridge, for support and hospitality during the programme \textit{Fractional Differential Equations (2022: FDE2)}
    where work on this paper was undertaken and presented in the
    framework of the Workshop "Optimal Control and Fractional Dynamics".

%%%%End of the main text
%%%%%%%%%%%%%%%

%%%%%%%APPENDICES
%%%%

        \end{document}